\documentclass[12pt,reqno]{amsart}
\usepackage{amssymb,latexsym}
\usepackage{verbatim}   
\theoremstyle{plain}
\newtheorem{theorem}{Theorem}

\newtheorem{lemma}[theorem]{Lemma}
\newtheorem{proposition}[theorem]{Proposition}

\theoremstyle{definition}

\theoremstyle{remark}

\begin{document}

\newcommand{\A}{{\mathfrak A}}
\newcommand{\B}{{\mathfrak B}}
\newcommand{\C}{{\mathfrak C}}
\newcommand{\Aut}{\textup{Aut}}
\newcommand{\BarHomo}{\,{}^{\overline{\hskip2.5mm}}\,}
\newcommand{\Char}{\textup{char}}
\newcommand{\diag}{\textup{diag}}
\newcommand{\E}{\mathbb{E}}
\newcommand{\EH}{\mathbb{E}H}
\newcommand{\Etimes}{\mathbb{E}^\times}
\newcommand{\End}{\textup{End}}
\newcommand{\F}{\mathbb{F}}
\newcommand{\FH}{\mathbb{F}H}
\newcommand{\Ftimes}{\mathbb{F}^\times}
\newcommand{\Frac}[2]{{\displaystyle\frac{#1}{#2}}}
\newcommand{\G}{\Gamma}
\newcommand{\Gal}{\textup{Gal}}
\newcommand\GL{\textup{GL}}
\newcommand\GLdE{\GL_d(\E)}
\newcommand\GLdF{\GL_d(\F)}
\newcommand\HatGamma{\widehat{\Gamma}}
\newcommand\HatPi{\widehat{\Pi}}
\newcommand\Hatpi{\widehat{\pi}}
\newcommand\Hom{\textup{Hom}}
\newcommand\im{\textup{im}}
\newcommand{\Inn}{\textup{Inn}}
\newcommand{\K}{\mathbb{K}}
\newcommand{\MatdE}{\E^{d\times d}}
\newcommand{\MatdF}{\F^{\kern1pt d\times d}}
\newcommand\Meataxe{{\sc Meataxe}}
\newcommand{\Q}{\mathbb{Q}}
\newcommand{\R}{\mathbb{R}}
\newcommand\tensor{\otimes}
\newcommand\Tr{\textup{Tr}}
\newcommand{\Z}{\mathbb{Z}}

\hyphenation{induced sub-fields Expo-sitions} 

\title[\tiny\upshape\rmfamily Writing representations over
proper sub-division rings]{}
\date{Draft printed on \today}

\begin{center}\large\sffamily\mdseries
 Writing representations over proper sub-division rings
\end{center}

\author{{\sffamily S.\,P. Glasby}}

\begin{abstract}
Let $\E$ be a division ring and $G$ a finite group of automorphisms of
$\E$ whose elements are distinct modulo inner automorphisms of
$\E$. Given a representation $\rho\colon\A\to\GLdE$ of an $\F$-algebra
$\A$, we give necessary and sufficient conditions for $\rho$ to be
{\it writable} over $\F=\E^G$, i.e. whether or not there exists a
matrix $A$ in $\GLdE$ that conjugates $\rho(\A)$ into $\GLdF$. We give an
algorithm for constructing an $A$, or proving that no $A$ exists. The
case of particular interest to us is when $\E$ is a field, and $\rho$
is absolutely irreducible. The algorithm relies on an explicit formula
for $A$, and a generalization of Hilbert's Theorem~90
(Theorem~3) that arises in Galois cohomology. The algorithm has
applications to the construction of absolutely irreducible group
representations (especially for solvable groups), and to the
recognition of one of the classes in Aschbacher's matrix group
classification scheme.
\end{abstract}

\maketitle
\centerline{\noindent Keywords: Hilbert's Theorem 90, proper sub-division ring}
\vskip2mm
\centerline{\noindent 2000 Mathematics subject classification: 20C40, 20C10}

\section{Introduction}\label{S:intro}

Throughout this paper $\E$ denotes a division ring, $G$ a finite group
of automorphisms of $\E$ whose elements are distinct modulo inner
automorphisms of $\E$, and $\F=\E^G$ is the sub-division ring fixed
elementwise by $G$. In the second half of this paper, we shall
specialize to the case when $\E:\F$ is a finite
Galois extension of {\it fields}. We view $\GLdE$ as the group of
invertible $d\times d$ matrices over $\E$. We say that a
representation $\rho\colon\A\to\GLdE$ of an $\F$-algebra
$\A$ {\it can be written over $\F$\/} if there exists an $A\in\GLdE$
such that
\[
 A^{-1}\rho(x)A\in\GLdF\qquad(x\in\A).
\]
The purpose of this paper is threefold: (1) to describe the connection
between Galois cohomology and the problem of writing $\rho$ over $\F$,
(2) to describe properties of a map $\Pi_C$ used to construct $A$,
and (3) to give an algorithm that takes as input an absolutely
irreducible $\rho$ and either constructs an $A$, or proves
that no such $A$ exists.

Section~2 describes briefly how $A$ gives rise to a certain function
$C\colon G\to\GLdE$ called 1-cocycle. The more interesting problem of
how $C$ gives rise to~$A$ is discussed in Section~3.  The heart of
this problem involves a generalization of Hilbert's Theorem 90: there
exists a matrix $A\in\GLdE$ such that $C_\alpha=A\alpha(A)^{-1}$ for
$\alpha\in G$. Equivalently, using the language of Galois cohomology,
it says that $H^1(G,\GLdE)=\{I\}$. This result was proved by Serre
\cite{S68} when $\E$ is a field, and by Nuss \cite{N97} when $\E$ is a
division ring. Neither the proof by Serre nor Nuss is constructive:
both proofs require modification in order to suggest an algorithm. We
shall give a completely elementary proof in Theorem~3 of these results
which suggests both a deterministic and a probabilistic algorithm for
constructing $A$. Although some of our results can be rephrased in terms of
Galois cohomology \cite{S68}, and descent theory for noncommutative
rings \cite{N97}, we prefer to state our results with minimal
background in terms of matrices over $\E$ and automorphisms.

Given a 1-cocycle $C\colon G\to\GLdE$, we can construct an
endomorphism $\Pi_C\colon\MatdE\to\MatdE$ of the algebra of
$d\times d$ matrices over $\E$. In Sections~3 and 4 we focus on
properties of $\Pi_C$. If $X$ is a random element of $\MatdE$, then
the probability that $A=\Pi_C(X)$ writes $\rho$ over $\F$ is at least
$\prod_{i=1}^\infty (1-2^{-i})>2/7$. After
Theorem~8 we shall assume that $\E$ is a (commutative) field. Different
choices for $X$ can give different choices for $A$, and a random $X$
can be a poor choice e.g. the entries of $A$ may be 100 digit
integers. We show in Theorem~10 that if $\E$ is a field and $|\F|\ge
d$, then we may take $X$ to be a scalar matrix. This
result, which is
best possible, appears to be helpful in producing ``nice'' conjugating
matrices $A$. Furthermore, whether $\lambda\in\E$ or $X\in\MatdE$, it
appears that the
probabilities $\textup{Prob}(\Pi_C(\lambda I)\textup{ invertible})$
and $\textup{Prob}(\Pi_C(X)\textup{ invertible})$ are very close.

Section~5 focuses on the case when $\rho$ is an absolutely
irreducible representation. In this case we construct a map $D\colon
G\to\GLdE$ and seek a function $\mu\colon G\to\Etimes$ such that $\mu
D$ is a 1-cocycle. The existence of $\mu$ determines whether or not
$\rho$ can be written over $\F$. If $\E$ is a cyclotomic number field,
then the existence of $\mu$ depends on the solutions to certain
equations in $\E$. We solve, if possible, certain norm equations, and
then solve equations in the group of units of the ring
of algebraic integers of $\E$.

Section~6 discusses some simple Las Vegas algorithms primarily for
computing $(q-1)$th roots, and solving norm equations in finite fields.
Section~7 gives examples arising from representations of
groups. Although our results apply to arbitrary $\F$-algebras $\A$,
the examples presented have $\A=\FH$ where $\FH$ is a group algebra
of a not necessarily finite group $H$. If $\sigma\colon H\to\GLdE$ is
a group representation, then $\sigma$ may be extended,
via a familiar argument, to a representation $\rho$ of the group
algebra $\A=\FH$. Of course, $\rho$ can be written over $\F$ precisely
when $\sigma$ can. The existence of a normal basis for $\E$ over $\F$
plays an important role in Section~7 and in Theorem~10.

Our work has been influenced by \cite{GH97}, which considers the case
when $G$ is cyclic, and by Br\"uckner's PhD thesis \cite{B98}.  In
\cite{B98} Br\"uckner independently discovers some results in
\cite{GH97}, and describes an unpublished result due to Plesken
\cite[Satz~3]{B98} which gives a necessary and sufficient condition
for an absolutely irreducible group representation over a field $\E$
to be writable over $\F$ where $\E:\F$ is a finite Galois extension of
fields. An algorithm is given in \cite[Lemma~7]{B98} for writing
$\rho$ over $\F$ when $G$ is cyclic. The proof contains
errors, however, all may be corrected. It involves choosing a random
column vector $x\in\E^{d\times 1}$ rather than choosing a random
matrix $X\in\MatdE$. This viewpoint motivated our Proposition~5.

In the sequel we will denote automorphisms of $\E$ by $\alpha$,
$\beta$, $\gamma$, elements of $\E$ by $\lambda$, $\mu$, $\nu$, and
representations of $\A$ by $\rho$, $\rho'$, $\sigma$.

\section{From $A$ to $C_\alpha$}

We shall say that {\it $\rho$ can be written over $\F$\/} if there
exists an $A\kern-1.6pt\in\kern-1.6pt\GLdE$ such that
\[
 A^{-1}\rho(x)A\in\GLdF\qquad(x\in\A).
\]
Our goal is to construct a conjugating matrix $A$, or prove that
one does not exist.

An automorphism $\alpha\in\Aut(\E)$ induces an automorphism,
also denoted $\alpha$, of the algebra $\MatdE$ of $d\times d$
matrices over $\E$: $(\mu_{i,j})\mapsto(\alpha(\mu_{i,j}))$.
Now $A$ writes $\rho$ over $\F$ if and only if
\[
 \alpha(A^{-1}\rho(x)A)=A^{-1}\rho(x)A\qquad(x\in\A,\alpha\in G).
\]
In subsequent equations, which hold for all $x\in\A$, we shall omit
the $x$'s and simply write
\[
 \alpha(A^{-1}\rho A)=A^{-1}\rho A\qquad(\alpha\in G).\tag{1}
\]
Therefore $C_\alpha:=A\alpha(A)^{-1}$ satisfies
\[
 C_\alpha^{-1}\rho\, C_\alpha=\alpha\circ\rho\qquad(\alpha\in G).\tag{2}
\]
Furthermore,
$A\alpha\beta(A)^{-1}=A\alpha(A)^{-1}\alpha(A\beta(A)^{-1})$ and so
\[
 C_{\alpha\beta}=C_\alpha\alpha(C_\beta)\qquad(\alpha,\beta\in G).\tag{3}
\]
We chose our automorphisms to act on the left, to avoid the
``twisted'' equation $C_{\alpha\beta}=C_\beta(C_\alpha)^\beta$,
which follows from $C_\alpha=A(A^\alpha)^{-1}$.

A map $C\colon G\to\GLdE$ defined by $\alpha\mapsto C_\alpha$
satisfying Eq.~(3) is called a {\it 1-cocycle}, and if there exists an
$A\in\GLdE$ such that $C_\alpha=A\alpha(A)^{-1}$ for all $\alpha\in
G$, then $C$ is called a {\it 1-coboundary}.
In summary, a necessary condition for $\rho$ to be writable over $\F$
is that there exist a 1-cocycle $C$ satisfying Eq.~(2). More
significantly, a 1-cocycle $C$ is a 1-coboundary, by a generalization
of Hilbert's Theorem~90, and there exist constructive methods for
finding $A$ from $C$, and hence for writing $\rho$ over~$\F$.

\section{From $C_\alpha$ to $A$}

The following result generalizes a well-known result of Artin
\cite[VIII\;\S4,\;Theorem~7]{L65} which
says that distinct characters $H\to\Etimes$ of a group $H$ with
values in a field $\E$, are linearly independent over $\E$.

\begin{lemma}
Let $\E$ be a division ring.
\begin{itemize}
\item[(a)] Let $\chi_1,\dots,\chi_n$ be group homomorphisms
$H\to\Etimes$ which are distinct modulo inner automorphisms of
$\E$. Then $\chi_1,\dots,\chi_n$ are linearly independent over $\E$.
\item[(b)] If $G$ is a finite subgroup of $\Aut(\E)$ whose elements
are distinct modulo $\Inn(\E)$, then the trace map
$\Tr\colon\E\to\F\colon\lambda\mapsto\sum_{\alpha\in G} \alpha(\lambda)$
is surjective.
\end{itemize}
\end{lemma}

\begin{proof}
(a) We shall view $\E$ as a left vector space over $\F$. The proof can
be modified for right $\F$-spaces.
Suppose that $\lambda_1\chi_1+\cdots+\lambda_n\chi_n=0$ where not
all $\lambda_i$ are zero, and $n$ is positive and minimal. Then
$n\ge2$ and each $\lambda_i$ is nonzero. If $h,k\in H$, then
\begin{align*}
 \lambda_1\chi_1(k)+\cdots+\lambda_n\chi_n(k)&=0,\\
 \lambda_1\chi_1(hk)+\cdots+\lambda_n\chi_n(hk)&=0.
\end{align*}
Premultiplying the first equation by
$\lambda_1\chi_1(h)\lambda_1^{-1}$ and subtracting the second
equation gives $\sum_{i=2}^n
\left(\lambda_1\chi_1(h)\lambda_1^{-1}\lambda_i-\lambda_i\chi_i(h)\right)\chi_i(k)=0$
for all $h,k\in H$. The minimality of $n$ implies that each coefficient
is zero. Therefore
$\chi_i(h)=\lambda_i^{-1}\lambda_1\chi_1(h)\lambda_1^{-1}\lambda_i$
for all $h\in H$, and $\chi_i$ is equivalent modulo $\Inn(\E)$
to $\chi_1$ for $i\ge2$, a contradiction.

(b) Let $\chi_1,\dots,\chi_n$ denote the elements of $G$ and let
$H=\Etimes$. By part~(a), $\chi_1,\dots,\chi_n$ are $\E$-linearly
independent and hence $\sum_{\alpha\in G} \alpha\ne 0$. Therefore the
$\F$-linear map $\Tr\colon\E\to\F$ is surjective.
\end{proof}

Assume we know matrices $C_\alpha\in\GLdE$ satisfying Eq.~(3). Theorem~3
shows how to construct $A\in\GLdE$ such that
$C_\alpha=A\alpha(A)^{-1}$ for $\alpha\in G$. It relies on the
following simple lemma.

\begin{lemma} Let $\E$ be a division ring, and let $G$ be a finite subgroup
of $\Aut(\E)$.
\begin{itemize}
\item[(a)] If $C_\alpha\in\MatdE$ satisfies
$C_{\alpha\beta}=C_\alpha+\alpha(C_\beta)$ for all $\alpha,\beta\in G$,
then $\Pi_C(X)=\sum_{\alpha\in G} C_\alpha+\alpha(X)$ satisfies
$C_\alpha+\alpha(\Pi_C(X))=\Pi_C(X)$ for all $X\in\MatdE$ and $\alpha\in G$.
\item[(b)] If $C_\alpha\in\GLdE$ satisfies
Eq.~(3), then $\Pi_C(X)=\sum_{\alpha\in G} C_\alpha \alpha(X)$ satisfies
$C_\alpha \alpha(\Pi_C(X))=\Pi_C(X)$ for all $X\in\MatdE$ and $\alpha\in G$.
\item[(c)] If $C_\alpha\in\GLdE$ satisfies Eq.~(3) and no two elements
of $G$ are equal modulo $\Inn(\E)$, then
then there exists a $\lambda\in\E$ such that the first column, $x$, of
$\Pi_C(I\lambda)$ is nonzero, and satisfies $C_\alpha \alpha(x)=x$ for all
$\alpha\in G$.
\end{itemize}
\end{lemma}

\begin{proof} We omit the proof of part (a) as it follows from the
proof of part (b) with products replaced by sums. It follows from
Eq.~(3) that
\[
 C_\alpha\alpha(\Pi_C(X))
   =C_\alpha\alpha\left(\sum_{\beta\in G}C_\beta \beta(X)\right)
   =\sum_{\alpha\in G}C_{\alpha\beta}\alpha\beta(X)
   =\Pi_C(X).
\]
Consider part (c).
Let $e$ be the column vector with 1 in the first
row, and zeroes elsewhere. Then $x=\Pi_C(I\lambda)e$, and by part (b)
\[
 C_\alpha\alpha(x)=C_\alpha \alpha(\Pi_C(I\lambda)e)
   =C_\alpha \alpha(\Pi_C(I\lambda))e
   =\Pi_C(I\lambda) e=x.
\]
Moreover, each of the column vectors of $C_\alpha \alpha(\lambda)$
are nonzero. By Lemma~1(b) the elements of $G$ are $\E$-linearly
independent. Hence there exists a $\lambda\in\E$ such that
$x=\sum_{\alpha\in G}  C_\alpha\alpha(\lambda) e\ne 0$.
\end{proof}

The sum $\sum C_\alpha \alpha(X)$ was considered in \cite{GH97}.
I have learned recently that this sum dates back to
Poincar\'e \cite[p. 159]{S68}. I attribute the following theorem to
Serre \cite[Prop.~3]{S68} when $\E$ is a field, and to Nuss
\cite[Theorem~B]{N97} when $\E$ is a division ring. We offer an
elementary proof conducive to practical implementation. A discussion
of non-matrix versions of Hilbert's Theorem~90 over division rings can
be found in \cite{L94}.

\begin{theorem} Let $\E$ be a division ring and $G$ a finite subgroup
of $\Aut(\E)$ whose elements are distinct modulo $\Inn(\E)$.
\begin{itemize}
\item[(a)] Let $C_\alpha\in\MatdE$, $\alpha\in G$. There exists an
$A\in\MatdE$ satisfying $C_\alpha=A-\alpha(A)$, $\alpha\in G$, if and
only if $C_{\alpha\beta}=C_\alpha+\alpha(C_\beta)$ for all
$\alpha,\beta\in G$.
\item[(b)] Let $C_\alpha\in\GLdE$, $\alpha\in G$. There exists an
$A\in\GLdE$ satisfying $C_\alpha=A\alpha(A)^{-1}$, $\alpha\in G$, if and
only if $C_{\alpha\beta}=C_\alpha\alpha(C_\beta)$ for all
$\alpha,\beta\in G$.
\end{itemize}
\end{theorem}

\begin{proof} The forward implication is straightforward for parts
(a) and~(b). The reverse implication follows from Lemma~2 for part
(a), and for part~(b) {\it provided} there exists and $X\in\MatdE$
such that $\Pi_C(X)$ is invertible. While it is clear that the
image of $\Pi_C$ contains {\it nonzero} matrices, it is more subtle that
$\im(\Pi_C)$ contains {\it invertible} matrices. We prove this second fact
via induction on $d$.

The result is true when $d=1$ by Lemma~2(c) since if $x\ne 0$, then
the $1\times 1$ matrix $[x]$ is invertible. Suppose that $d>1$ and
that the result is true for dimension $d-1$. By Lemma~2(c) there
exists an invertible matrix $Y$ with first column $x$, satisfying
$C_\alpha \alpha(x)=x$ for all $\alpha\in G$. Therefore,
\[
 Y^{-1}C_\alpha \alpha(Y)
  =\begin{pmatrix} 1&y_\alpha\\0&C'_\alpha\end{pmatrix}\qquad (\alpha\in G)
\]
where $C'_\alpha\in\GL_{d-1}(\E)$. Since $Y^{-1}C_\alpha \alpha(Y)$
satisfies Eq.~(3), so too does $C'_\alpha$. By induction, there exists
an $A'\in\GL_{d-1}(\E)$ satisfying $C'_\alpha \alpha(A')=A'$ for all
$\alpha\in G$. Thus
\[
 \begin{pmatrix} 1&0\\0&A'\end{pmatrix}^{-1}
 Y^{-1}C_\alpha \alpha(Y)
 \alpha\begin{pmatrix} 1&0\\0&A'\end{pmatrix}
 =\begin{pmatrix} 1&z_\alpha\\0&I\end{pmatrix}=:C''_\alpha\qquad(\alpha\in G).
\]
Since $C''_\alpha$ satisfies Eq.~(3), the $z_\alpha$ satisfy
$z_{\alpha\beta}=z_\alpha+\alpha(z_\beta)$ for all $\alpha,\beta\in
G$. By part (a) there exists a $1\times(d-1)$ vector $w$ such that
$z_\alpha=w-\alpha(w)$ for all $\alpha\in G$. Therefore,
$A=Y\begin{pmatrix} 1&0\\0&A'\end{pmatrix}
   \begin{pmatrix} 1&w\\0&I\end{pmatrix}$.
\end{proof}

Lemma~2(b) entreats us to study the maps
$\Pi_C, \Gamma_\alpha\colon\MatdE\to\MatdE$ defined by
\[
 \Pi_C(X)=\sum_{\alpha\in G} C_\alpha \alpha(X)\quad\text{and}\quad
 \Gamma_\alpha(X)=C_\alpha \alpha(X)-X.
\]
When $\Char(\E)\nmid |G|$, it is convenient to also define $\pi_C$ by
$\pi_C=|G|^{-1}\Pi_C$.
The matrix $A$ in Theorem~3 satisfying $C_\alpha=A\alpha(A)^{-1}$ is
far from unique. Indeed the matrix $AY$, where $Y\in\GLdF$, has the same
property. It is useful to regard $\MatdE$ as a right $\MatdF$-module,
where the scalar action is right matrix multiplication.

\begin{proposition}
Let $C\colon G\to\GLdE$ be a 1-cocycle where $\E$ is a division ring
and $G$ is a finite subgroup of $\Aut(\E)$.
\begin{itemize}
\item[(a)] The maps $\Pi_C$ and $\Gamma_\alpha$ are right
$\MatdF$-homomorphisms satisfying
$\Pi_C\circ\Gamma_\alpha=\Gamma_\alpha\circ\Pi_C=0$ and $\Pi_C^2=|G|\Pi_C$.
\item[(b)] If $\Char(\E)\nmid |G|$, then $\pi_C^2=\pi_C$ and so
$\MatdE=\im(\pi_C)\dotplus\ker(\pi_C)$ where
$\ker(\pi_C)=\im(1-\pi_C)$. Moreover, if $\pi_C(X)=XY$ where
$Y\in\GLdF$, then $\pi_C(X)=X$.
\item[(c)] If $C_\alpha=A\alpha(A)^{-1}$ for all $\alpha\in G$, then
$\Pi_C(X)=A\Tr(A^{-1}X)$ where $\Tr\colon\MatdE\kern-2pt\to\MatdF$ is the
trace function: $X\mapsto\sum_{\alpha\in G}\alpha(X)$. Moreover,
$\Pi_C(A\lambda)=A\Tr(\lambda)$, $\Pi_C(A)=|G|A$ and $\pi_C(A)=A$.
\item[(d)] Let $Y\in\GLdE$ be fixed, and let $D\colon G\to\GLdE$ be
defined by $D_\alpha=Y^{-1}C_\alpha \alpha(Y)$. Then $D_\alpha$
satisfies \textup{Eq.~(3)}, and
\[
 \Pi_D(X)=Y^{-1}\Pi_C(YX).
\]
\end{itemize}
\end{proposition}

\begin{proof}
(a) It is clear that $\Pi_C(X_1+X_2)=\Pi_C(X_1)+\Pi_C(X_2)$
and $\Pi_C(XY)=\Pi_C(X)Y$ for all $Y\in\MatdF$. Thus $\Pi_C$, and
similarly $\Gamma_\alpha$, are right $\MatdF$-module homomorphisms.
Lemma~2(b) shows that $\Gamma_\alpha\circ\Pi_C=0$, and the
following argument shows that $\Pi_C\circ\Gamma_\beta=0$:
\[
 \Pi_C(C_\beta \beta(X))
 =\sum_{\alpha\in G} C_\alpha \alpha(C_\beta \beta(X))
 =\sum_{\alpha\in G} C_{\alpha\beta} \alpha\beta(X)=\Pi_C(X).
\]
In addition, by the above equation:
\[
 \Pi_C^2(X)=\sum_{\beta\in G} \Pi_C(C_\beta \beta(X))
 = \sum_{\beta\in G} \Pi_C(X)=|G|\Pi_C(X).
\]

(b) Multiplying the equation $\Pi_C^2=|G|\Pi_C$ by $|G|^{-2}$ gives
$\pi_C^2=\pi_C$. Standard arguments show that
$\MatdE=\im(\pi_C)\dotplus\ker(\pi_C)$. If $\pi_C(X)$ equals $XY$,
then by part~(a)
\[
 XY=\pi_C(X)=\pi_C^2(X)=\pi_C(XY)=\pi_C(X)Y=XY^2.
\]
Postmultiplying by $Y^{-1}$ gives $X=XY$. Thus $\pi_C(X)=X$.

Consider part~(c):
\[
 \Pi_C(X)=\sum_{\alpha\in G} A\alpha(A)^{-1}\alpha(X)
 =A\sum_{\alpha\in G} \alpha(A^{-1}X)=A\Tr(A^{-1}X).
\]
Setting $X=A\lambda$ shows $\Pi_C(A\lambda)=A\Tr(I\lambda)=A\Tr(\lambda)$, and
setting $\lambda=1$ shows $\Pi_C(A)=|G|A$ and $\pi_C(A)=A$.
Part~(d) is straightforward. (The 1-cocycles $C$ and $D$ are called
{\it cohomologous}.)
\end{proof}

The endomorphisms $\Pi_C,\Gamma_\alpha$ of $\MatdE$ give rise to
endomorphisms $\HatPi_C$, $\HatGamma_\alpha$ of the space $\E^{d\times1}$
of $d\times1$ column vectors:
\[
 \HatPi_C(x)=\sum_{\alpha\in G} C_\alpha \alpha(x),\quad
 \HatGamma_\alpha(x)=C_\alpha \alpha(x)-x\qquad(\alpha\in G).
\]
When $\Char(\E)\nmid |G|$, it is convenient to also define $\Hatpi_C$ by
$\Hatpi_C=|G|^{-1}\HatPi_C$.
If $x\in\E^{d\times1}$ is the first column of $X\in\MatdE$,
and $Y=\textup{diag}(1,0,\dots,0)$, then the first columns of
$\Pi_C(XY)=\Pi_C(X)Y$ and $\Gamma_\alpha(XY)=\Gamma_\alpha(X)Y$
are $\HatPi_C(x)$ and $\HatGamma_\alpha(x)$ respectively.

It is worth recording some simple generalizations of
Prop.~4(a,b,c) such as:
$\HatGamma_\alpha\circ\HatPi_C=\HatPi_C\circ\HatGamma_\alpha=0$,
$\HatPi_C^2=|G|\HatPi_C$,
$\E^{d\times1}=\im(\Hatpi_C)\dotplus\ker(\Hatpi_C)$ and
$\HatPi_C(x)=A\Tr(A^{-1}x)$ where $\Tr$ denotes the trace map
$\E^{d\times1}\to\F^{d\times1}$.

\begin{proposition}
Let $C\colon G\to\GLdE$ be a 1-cocycle where $\E$ is a division ring
and $G$ is a finite subgroup of $\Aut(\E)$ whose elements are
distinct modulo $\Inn(\E)$. Let $S$ be a generating set for $G$, and let
$\F=\E^G$.
\begin{itemize}
\item[(a)] $\im(\HatPi_C)=\bigcap_{\alpha\in S} \ker(\HatGamma_C)$ is the
$\F$-linear span of the columns of any matrix $A$ satisfying
$C_\alpha=A\alpha(A)^{-1}$ for all $\alpha\in G$.
\item[(b)] If $\Char(\E)\nmid |G|$, then
$\ker(\HatPi_C)=\sum_{\alpha\in S} \im(\HatGamma_\alpha)$.
\item[(c)] If $\alpha\ne1$, then $\im(\HatGamma_\alpha)$ spans
$\E^{d\times1}$ as an $\E$-space.
\item[(d)] If $0\ne x\in\ker(\HatPi_C)$, then
$x\E\not\subseteq\ker(\HatPi_C)$.
\end{itemize}
\end{proposition}

\begin{proof}
(a) $\HatGamma_\alpha\circ\HatPi_C=0$, implies
$\im(\HatPi_C)\subseteq\bigcap_{\alpha\in
S}\ker(\HatGamma_\alpha)$. Conversely, if $x\in\bigcap_{\alpha\in
S}\ker(\HatGamma_\alpha)$, then $C_\alpha\alpha(x)=x$ for $\alpha\in
S$. It follows from Eq.~(3) that $C_\alpha\alpha(x)=x$ for $\alpha\in
G$. Thus
\[
 \Pi_C(x\lambda)=\sum_{\alpha\in G} C_\alpha\alpha(x)\alpha(\lambda)
   =\sum_{\alpha\in G} x\alpha(\lambda)=x\Tr(\lambda).
\]
By Lemma~1(b), there exists a $\lambda\in\E$ such that
$\Tr(\lambda)=1$. Thus $x\in\im(\HatPi_C)$ and so
$\im(\HatPi_C)=\bigcap_{\alpha\in S}\ker(\HatGamma_\alpha)$. It
follows from Prop.~4(c) that $\im(\HatPi_C)=A\MatdF$, and so
$\im(\HatPi_C)$ is the $\F$-linear span of columns of $A$.

(b) $\HatPi_C\circ\HatGamma_\alpha=0$, implies
$\ker(\HatPi_C)\supseteq\sum_{\alpha\in S} \im(\HatGamma_\alpha)$.
It follows from Eq.~(3) that
\[
 C_{\alpha\beta}\alpha\beta(x)-x
   =[C_{\alpha}\alpha(C_\beta\beta(x))-C_{\beta}\beta(x)]+
   [C_{\beta}\beta(x)-x].
\]
Hence $\im(\HatGamma_{\alpha\beta})\subseteq
\im(\HatGamma_{\alpha})+\im(\HatGamma_{\beta})$ and
$\sum_{\alpha\in G} \im(\HatGamma_\alpha)=
\sum_{\alpha\in S} \im(\HatGamma_\alpha)$.
Conversely, if $x\in\ker(\HatPi_C)$, then
$\sum_{\alpha\in G}C_\alpha \alpha(x)=0$ and hence
\[
 x=\Tr(|G|^{-1}x)=\sum_{\alpha\in G}\HatGamma_{\alpha}(|G|^{-1}x)\in
 \sum_{\alpha\in G} \im(\HatGamma_\alpha)=
 \sum_{\alpha\in S} \im(\HatGamma_\alpha).
\]
Thus $\ker(\HatPi_C)=\sum_{\alpha\in S} \im(\HatGamma_\alpha)$ as
desired.

(c) Suppose that Let $\phi\colon\E^{d\times1}\to\E$ be an $\E$-linear map
containing $\im(\HatGamma_\alpha)$ in its kernel. Then for all
$x\in\E^{d\times1}$ and $\lambda\in\E$:
\[
 0=\phi(\HatGamma_\alpha(x\lambda))=\phi(C_\alpha
 \alpha(x))\alpha(\lambda) - \phi(x)\lambda.
\]
Since $\alpha\ne 1$ it follows from Lemma~1(a) that $\phi(x)=0$ for
all $x$ and hence $\phi=0$. This proves that the $\E$-linear span of
$\im(\HatGamma_\alpha)$ equals $\E^{d\times1}$, and hence
$\dim_\F(\im(\HatGamma_\alpha))\ge d$.

(d) Suppose that $0\ne x\in\ker(\HatPi_C)$. If $\HatPi_C(x\lambda)=0$
for all $\lambda\in\E$, then $\sum_{\alpha\in G}C_\alpha
\alpha(x)\alpha(\lambda)=0$. Since $C_\alpha \alpha(x)\ne 0$, this
contradicts Lemma~1(a). Thus $x\E\not\subseteq\ker(\HatPi_C)$ as claimed.
\end{proof}

In light of Prop.~5(a) the assumption in Prop.~5(b)
that $\Char(\E)\nmid|G|$ may be unnecessary.

\begin{proposition} Let $(\lambda_\alpha)_{\alpha\in G}$ be an
$\F$-basis for $\E$, and let $E_{i,j}\in\MatdE$ be the matrix with 1
in the $(i,j)$th entry and zeroes elsewhere. Then $\MatdE$ is a freely
generated as a right $\MatdF$-module by $E_{i,1}\lambda_\alpha$,
$\alpha\in G$, $i=1,\dots,d$.
\end{proposition}

\begin{proof}
By taking $\F$-linear combinations of $E_{i,1}\lambda_\alpha$ gives a
matrix with arbitrary first column. Taking $\MatdF$-multiples gives
every element of $\MatdE$. The fact that the $E_{i,1}\lambda_\alpha$
freely generate $\MatdE$ follows from the observation that
$E_{i,1}\MatdF$ comprises matrices with all rows zero except the $i$th,
and the $i$th row can be an arbitrary vector in $\F^{\,1\times d}$.
\end{proof}

It follows from Theorem~3 and the above proposition that an invertible
matrix can be found by taking $\MatdF$-linear combinations of the
matrices $\Pi_C(\lambda_\alpha E_{i,1})$. Since each
$\Pi_C(\lambda_\alpha E_{i,1})$ is
singular (unless $d=1$), it is better to consider $\MatdF$-linear
combinations of $\Pi_C(\lambda_\alpha D^i)$ where $D$ is the permutation
matrix corresponding to the $d$-cycle $(1,2,\dots,d)$. A simple
argument shows that the $\lambda_\alpha D^i$ generate $\MatdE$ as a
$\MatdF$-module,
although not freely. In practice $\MatdF$-linear combinations are not
necessary as $\Pi_C(\lambda_\alpha D^i)$ is commonly
invertible. Thus we typically do not evaluate $\Pi_C(X)$ at a random
matrix $X$. Doing so can result in ``bad'' matrices $A=\Pi_C(X)$,
e.g. with 100 digit integer entries. More significantly, the matrices
$A^{-1}\rho(x)A$ can be ``bad''. Choosing $X$ to be a scalar matrix
seems to result in ``good'' matrices $\Pi_C(X)$. This imprecise
statement has some theoretical underpinning in Theorem~10.

\section{Invertible elements in $\textup{im}(\Pi_c)$}

The primary aim of this section is to prove in Theorem~10 that if
$|\F|\ge d$ there exists a $\lambda\in\E$ such that $\Pi_C(I\lambda)$
is invertible. We show in Theorem~~8 that the assumption
$|\F|\ge d$ is best possible by considering a special case when
$A$, and hence each $C_\alpha$, is upper-triangular.

We need a preliminary lemma.

\begin{lemma} Let $V$ be a vector space over a division ring $\F$.
If $V$ is a union of $m$ proper subspaces, then $\dim_\F(V)\ge 2$ and
$|\F|< m$. Conversely, if $\dim_\F(V)\ge 2$ and $\F$ is finite, then
$V$ is a union of $|\F|+1$ proper subspaces.
\end{lemma}

\begin{proof}
The proof in \cite[Problem~24]{H95} generalizes to division rings. If
$\dim_\F(V)\ge 2$, then $V=H_\infty\cup\bigcup_{\lambda\in\F}
H_\lambda$ where $H_\infty$ is the hyperplane $x_1=0$ and $H_\lambda$
the hyperplane $\lambda x_1+x_2=0$. Thus $V$ is a union of $|\F|+1$
proper subspaces.
\end{proof}

\begin{theorem} Let $C\colon G\to\GLdE$ be a 1-cocycle where $G$ and
$\E$ are as in Theorem~3. Suppose that $C_\alpha=A\alpha(A)^{-1}$,
$\alpha\in G$, where $A\in\GLdE$ is upper-triangular and $\F=\E^G$. If
$|\F|\ge d$, then there exists a $\lambda\in\E$ such that
$\Pi_C(I\lambda)$ is invertible. Moreover, if $|\F|<d$,
then there exists an upper-triangular matrix $A\in\GLdE$ such that
$\Pi_C(I\lambda)$ is singular for all $\lambda\in\E$.
\end{theorem}

\begin{proof} It follows from Prop.~4(c) that $\Pi_C(I\lambda)$
is invertible if and only if $\Tr(A^{-1}\lambda)$ is invertible. If
$a_{i,i}$ denotes the $(i,i)$th entry of $A$, then
$\Tr(A^{-1}\lambda)$
is upper-triangular with $(i,i)$th entry $\Tr(a_{i,i}^{-1}\lambda)$.
Let $K(a_{i,i}^{-1})$ denote the kernel of the map
$\lambda\mapsto\Tr(a_{i,i}^{-1}\lambda)$. By Lemma~1(b) the $\F$-subspace
$K(a_{i,i}^{-1})$ of $\E$ has codimension~1. If $|\F|\ge d$, then $\E$ is not a
union of $d$ proper subspaces by Lemma~7. Thus there exists a
$\lambda\in\E$ not in
$\bigcup_{i=1}^d K(a_{i,i}^{-1})$. Since $\Tr(a_{i,i}^{-1}\lambda)\ne 0$ for
each $i$, it follows that $\Pi_C(I\lambda)$ is invertible.

Conversely, suppose that $|\F|<d$. Then $\E$ is a union of $|\F|+1$
proper subspaces, so we may choose
$a_{1,1}^{-1},\dots,a_{d,d}^{-1}\in\Etimes$ such that
$\E=\bigcup_{i=1}^d K(a_{i,i}^{-1})$. Then for each $\lambda\in\E$
at least one diagonal entry of the upper-triangular matrix
$\Tr(A^{-1}\lambda)$ is zero. Put differently,
$\Pi_C(I\lambda)$ is singular for all $\lambda\in\E$.
\end{proof}

\noindent {\bf Assumption:} We shall henceforth assume that $\E$ is a
{\it field}.

Theorem~10 generalizes Theorem~8 to deal with arbitrary $d\times
d$ matrices $A$. Its proof assumes that $\E$ is a field, and depends on
the following well-known result.

\begin{lemma} Let $f$ be an element of the polynomial ring $\F[x_1,\dots,x_n]$
such that $f(a_1,\dots,a_n)=0$ for all $(a_1,\dots,a_n)\in\F^n$.
\begin{itemize}
\item[(a)] If the degree of $f$ in each variable is less than $|\F|$,
then $f=0$.
\item[(b)] If the degree of $f$ is at most $q$ where
$|\F|=q$, then there exists $\nu_1,\dots,\nu_n\in\F$
such that $f(x_1,\dots,x_n)=\sum_{i=1}^n \nu_i(x_i^q-x_i)$.
\end{itemize}
\end{lemma}

\begin{proof}
(a) See \cite[Chapter V, Theorem~5]{L65} (and \cite[Corollary~3]{L65})
for the case when $\F$ is finite (and $\F$ is infinite).
Consider part (b). Recall that the degree of a nonzero polynomial is the
maximum degree of a monomial summand, and $\deg(x_1^{k_1}\cdots
x_n^{k_n})=k_1+\cdots+k_n$. The result is true
when $n=1$. Suppose that $n>1$ and $f=\sum_{i=0}^q f_ix_n^{q-i}$ where
$f_i$ is a polynomial in $x_1,\dots,x_{n-1}$ of degree at most $i$.
Fix $(a_1,\dots,a_{n-1})\in\F^{n-1}$ and consider
$f(a_1,\dots,a_{n-1},x^n)$. By the $n=1$ case,
$f_i(a_1,\dots,a_{n-1})=0$ for $i=1,\dots,q-2$ and
$f_0=-f_{q-1}=\nu_n$ is a constant polynomial. By part~(a), $f_i=0$ for
$i=1,\dots,q-2$ and by induction there exist
$\nu_1,\dots,\nu_{n-1}\in\F$ such that $f_q=\sum_{i=1}^{n-1}
\nu_i(x_i^q-x_i)$. In summary, $f=\sum_{i=1}^n \nu_i(x_i^q-x_i)$.
\end{proof}

The reader may like to compare Lemma~9(b) with a theorem due to
Chevalley \cite[\S1.7, Theorem~2]{S71}.

\begin{theorem}
Let $\E$ be a field, and $\E:\F$ a finite Galois extension with group
$G$. Suppose that $C\colon G\to\GLdE$ is a 1-cocycle and $|\F|\ge
d$. Then there exists a $\lambda\in\E$ such that $\Pi_C(I\lambda)
=\sum_{\alpha\in G} C_\alpha \alpha(\lambda)$ is invertible.
\end{theorem}

\begin{proof}
By Theorem~3 there exists an invertible matrix $A$ satisfying
$C_\alpha=A\alpha(A)^{-1}$, $\alpha\in G$. By Prop.~4(c),
$\Pi_C(\lambda I)=A\Tr(\lambda A^{-1})$. Thus $\Pi_C(\lambda I)$ is invertible
precisely when $\Tr(\lambda A^{-1})$ is invertible. Our problem can be
rephrased: Given $X\in\GLdE$, find $\lambda\in\E$ such that
$\Tr(\lambda X)$ is invertible.

By \cite[Theorems~7.4.2, 8.7.2]{R95} there exists $\zeta\in\E$ such
that $(\alpha(\zeta))_{\alpha\in G}$ is a basis for $\E$ over
$\F$ (such a basis is called a {\it normal basis}).
Now $\Tr(\zeta)\in\Ftimes$ by Lemma~1(b). By replacing $\zeta$ by
$\Tr(\zeta)^{-1}\zeta$ we
may additionally assume that $\Tr(\zeta)=1$. A typical element of $\E$
has the form $\sum_{\alpha\in G} x_\alpha \alpha(\zeta)$ where
$x_\alpha\in\F$. Write
\[
 x_{i,j}=\sum_{\alpha\in G}
 x_\alpha^{i,j}\alpha(\zeta)\quad\text{and}\quad
 \lambda=\sum_{\beta\in G}\lambda_\beta \beta(\zeta)
\]
where $x_{i,j}$ denotes the $(i,j)$th entry of $X$. We shall view the
$x_\alpha^{i,j}$ as {\it elements} of $\F$, and the $\lambda_\beta$ as
algebraically independent commuting {\it variables} that are fixed by $G$.

Let $(\mu_{\alpha,\beta})$ be the matrix of the $\F$-linear
transformation $\E\to\E$ defined by $\lambda\mapsto\zeta\lambda$. That is,
\[
 \zeta\alpha(\zeta)=\sum_{\beta\in G}\mu_{\alpha,\beta}\beta(\zeta)
 \quad\quad (\mu_{\alpha,\beta}\in\F).\tag{4}
\]
Then
\begin{align*}
 x\lambda&=\left(\sum_{\alpha}x_\alpha \alpha(\zeta)\right)
 \left(\sum_{\beta}\lambda_\beta \beta(\zeta)\right)
 =\sum_{\alpha,\beta}
  x_\alpha\lambda_\beta \alpha(\zeta \alpha^{-1}\beta(\zeta))\\
 &=\sum_{\alpha,\beta,\gamma}
  x_\alpha\lambda_\beta\mu_{\alpha^{-1}\beta,\gamma}\; \alpha\gamma(\zeta).
\end{align*}
Replacing $\alpha\gamma$ by $\gamma$ gives
$x\lambda=\sum x_\alpha\lambda_\beta
\mu_{\alpha^{-1}\beta,\alpha^{-1}\gamma}\;\gamma(\zeta)$.
Our normalization implies that $\Tr(\gamma(\zeta))=1$, and
hence
\[
 \Tr(x\lambda)=\sum_{\alpha}\left(\sum_{\beta,\gamma}
  \mu_{\alpha^{-1}\beta,\alpha^{-1}\gamma}\lambda_\beta\right)x_\alpha.
 \qquad\tag{5}
\]
Abbreviate the above inner sum by $z_\alpha$. Then
\begin{align*}
 z_\alpha&=\sum_{\beta}\left(
   \sum_{\gamma}\mu_{\alpha^{-1}\beta,\alpha^{-1}\gamma}
   \right)\lambda_\beta
   =\sum_{\beta}\Tr(\zeta \alpha^{-1}\beta(\zeta))\lambda_\beta\\
   &=\sum_{\beta}\Tr(\alpha(\zeta) \beta(\zeta))\lambda_\beta.
   \qquad\tag{6}
\end{align*}
Replacing $x_\alpha$ in Eq.~(5) by $x^{i,j}_\alpha$ shows
\[
  \det\Tr(X\lambda)=\det(x_{i,j}\lambda)
  =\det\left(\sum_\alpha z_\alpha x^{i,j}_\alpha\right).
\]
This determinant is a polynomial in the variables $z_\alpha$ which is
either the zero polynomial, or is homogeneous of degree~$d$ in the
$z_\alpha$. Specifically,
\[
 \det\left(\sum_\alpha z_\alpha x^{i,j}_\alpha\right)
 =\sum p_{\{\alpha_1,\dots,\alpha_d\}}
   z_{\alpha_1}\cdots z_{\alpha_d}\qquad\tag{7}
\]
where the sum is taken over all orbits of the symmetric group $S_d$ on
the group $G^d$. Such orbits are in bijective correspondence with the
multisets $\{\alpha_1,\dots,\alpha_d\}$ of $G$ having at least one,
and at most $d$, distinct elements. We view the coefficient
$p_{\{\alpha_1,\dots,\alpha_d\}}$ of $z_{\alpha_1}\cdots
z_{\alpha_d}$ as an element of $\F$, not a polynomial over $\F$ in the
$x^{i,j}_\alpha$.

The matrix $(\Tr(\alpha(\zeta)\beta(\zeta)))_{\alpha,\beta\in G}$ is
invertible (see \cite[\S7.2]{R95}), and its determinant equals the
discriminant $\prod_{\alpha\ne\beta} (\alpha(\zeta)-\beta(\zeta))$ of the
minimal polynomial $\prod_\alpha (t-\alpha(\zeta))$ of $\zeta$ over
$\F$. By Eq.~(6) as $(\lambda_\beta)$ runs through the
vectors in the vector space $\F^{|G|}$, $(z_\alpha)$ does the same.

The determinant $\det(X)=\det(\sum_\alpha x^{i,j}_\alpha\alpha(\zeta))$
can be evaluated using the same reasoning used for Eq.~(7).
Replacing $z_\alpha$ by $\alpha(\zeta)$ in Eq.~(7) shows
\[
 \det(X)=\sum p_{\{\alpha_1,\dots,\alpha_d\}}
  \alpha_1(\zeta)\cdots\alpha_d(\zeta).\qquad\tag{8}
\]

Let us assume that $X$ is fixed and that $\det\Tr(X\lambda)=0$ for
all $\lambda\in\E$ (or equivalently, all $(\lambda_\beta)\in\F^{|G|}$).
By virtue of the previous paragraph, this says that the polynomial
Eq.~(7) is zero for all $(z_\alpha)\in\F^{|G|}$. If $|\F|>d$, then
Lemma~9(a) implies that each $p_{\{\alpha_1,\dots,\alpha_d\}}$ equals zero. By
Eq.~(8), $\det(X)=0$. In summary, we have proved that if $|\F|>d$ and
$\det(X)\ne 0$, then there exists a $\lambda\in\E$ such that $\det\Tr(\lambda
X)\ne 0$.

Finally, suppose that $|\F|=d$ is finite and $\det\Tr(X\lambda)=0$ for
all $\lambda\in\E$. By Lemma~9(b), $\det\Tr(X\lambda)=\sum_{\alpha\in G}
\nu_\alpha(z_\alpha^{|\F|}-z_\alpha)$. Since
this polynomial is not homogeneous, each $\nu_\alpha$
is zero. Thus each $p_{\{\alpha_1,\dots,\alpha_d\}}$ equals zero, and
$\det(X)=0$ by Eq.~(8). This completes the proof.
\end{proof}

In the light of Theorem~8, one may suspect that Theorem~10 holds more
generally: namely when $\E$ is a division ring.

\section{Algorithmic considerations}

Henceforth assume that $\rho\colon\A\to\GLdE$ is an {\it
absolutely irreducible} representation, and $\langle S\mid R\rangle$ is
a finite presentation of $G$.

If $\rho$ can be written over $\F$, then there exist matrices
$D_\alpha\in\GLdE$ satisfying
\[
 D_\alpha^{-1}\rho D_\alpha=\alpha\circ\rho\qquad\qquad(\alpha\in G).\tag{9}
\]
There are a variety of methods for calculating the $D_\alpha$, or proving
that some do not exist. These include (a) using the \Meataxe\ algorithm
\cite{HR94,P98}, (b)~solving $d^2|G|$ homogeneous linear equations
over $\F$ in $d^2|G|$ unknowns, and (c) averaging over a chain
$\A=\A_1\supset\cdots\supset\A_{n+1}=\{0\}$ of $\F$-algebras
where the indices $|\A_i:\A_{i+1}|$ are ``small''. 

If $\rho$ can be written over $\F$, then by absolute irreducibility there
exists a function $\mu\colon G\to\Etimes$ such that $C=\mu D$ is a
1-cocycle. It suffices to know $C_\alpha$ for $\alpha\in S$, because
Eq.~(3) allows us to compute $C_\gamma$ for $\gamma\in G$.
Suppose henceforth that we have computed matrices
$D_\alpha$, $\alpha\in S$, that satisfy Eq.~(9).
Now the $C_\alpha$, $\alpha\in S$, satisfy the relations $R$ for $G$, and
in general the $D_\alpha$ will not. The relations give rise to a
system of $|R|$ equations that the scalars $\mu_\alpha$ must satisfy. If these
equations can not be solved in $\Etimes$, then $\rho$ can not be
written over $\F$, otherwise it can by Section~2.  We shall say more
about the equations that the $\mu_\alpha$ satisfy.

Two important applications of this work are (a) when $\E$ is a subfield
of a cyclotomic field, and (b) when $\E$ is a finite field. In these
cases $G$ is abelian, or cyclic and we assume that
$G$ has a presentation:
\[
  G=\langle\alpha_1,\dots,\alpha_s\mid \alpha_i^{m_i}=1,1\le i\le s,
       [\alpha_j,\alpha_i]=1, 1\le i<j\le s\rangle
\]
where $[\alpha,\beta]$ denotes the commutator
$\alpha^{-1}\beta^{-1}\alpha\beta$. We shall not necessarily assume
that $m_1|m_2|\cdots|m_s$.

The power relations and the commutator relations give different
equations that the $\mu_\alpha$ must satisfy. Consider first relations
of the form $\alpha^m=1$. It follows from Eq.~(3) that
$C_{\alpha^m}=C_\alpha \alpha(C_\alpha)\cdots\alpha^{m-1}(C_\alpha)$
and hence that
\begin{align*}
 D_\alpha\alpha(D_\alpha)\cdots\alpha^{m-1}(D_\alpha)&=\lambda_{\alpha}I
 \quad\text{where}\\
 \mu_\alpha\alpha(\mu_\alpha)\cdots\alpha^{m-1}(\mu_\alpha)
 &=\lambda_{\alpha}^{-1}\qquad(\alpha\in S).\tag{10}
\end{align*}
Given a subgroup $A$ of $G$, define the norm map
$N_A\colon\Etimes\to(\E^A)^\times$ by
$N_A(\lambda)=\prod_{\alpha\in A} \alpha(\lambda)$. Then Eq.~(10) says:
$N_{\langle\alpha\rangle}(D_\alpha)=\lambda_{\alpha}I$ where
$N_{\langle\alpha\rangle}(\mu_\alpha)=\lambda_{\alpha}^{-1}$
for some $\mu_\alpha\in\Etimes$.
A necessary condition is that
$\alpha(\lambda_\alpha)=\lambda_\alpha$.
When $\E$ is finite, $N_{\langle\alpha\rangle}$ is surjective, and
this necessary condition is sufficient to guarantee a solution for
$\mu_\alpha$. By contrast, when $\E$ is infinite the equation
$N_{\langle\alpha\rangle}(\mu_\alpha)=\lambda_\alpha^{-1}$ may have no solution
(c.f. Section~7, Example~1). There are a variety of algorithms for
solving for $\mu_\alpha$ when $\E$ is finite, see for example
Section~6. Different algorithms are required in the the case when $\E$
is a number field, see for example \cite{F97} and \cite{S02}. Assume
henceforth that the equations
$N_{\langle\alpha\rangle}(\mu_\alpha)=\lambda_\alpha^{-1}$ can be solved. By
replacing $D_\alpha$ by $\mu_\alpha^{-1}D_\alpha$ we will henceforth
assume that $\lambda_\alpha=1$ for $\alpha\in S$. We shall now seek a
function $\nu$ such that $\nu D$ is a 1-cocycle.

Consider now equations arising from commutators
$[\alpha,\beta]=1$ where $\alpha,\beta\in S$. Applying Eq.~(3) twice
gives
\[
 C_\alpha\alpha(C_\beta)=C_{\alpha\beta}
 =C_{\beta\alpha}=C_\beta\beta(C_\alpha)\qquad(\alpha,\beta\in G).
\]
Substituting $C_\alpha=\nu_\alpha D_\alpha$ into
$\alpha(C_\beta)^{-1}C_\alpha^{-1}C_\beta\beta(C_\alpha)=I$ gives
\begin{align*}
 \alpha(D_\beta)^{-1}D_\alpha^{-1}D_\beta\beta(D_\alpha)
 &=\lambda_{\alpha,\beta}I\quad\text{where}\\
 \beta(\nu_\alpha)^{-1}\nu_\beta^{-1}\nu_\alpha\alpha(\nu_\beta)
 &=\lambda_{\alpha,\beta}\qquad(\alpha,\beta\in S).
 \quad\tag{11}
\end{align*}
Let $K_A$ and $I_A$ denote the kernel and image of the norm map $N_A$.
As we are assuming that $\lambda_\alpha=\lambda_\beta=1$
it follows that $\nu_\alpha\in K_A$ and $\nu_\beta\in K_B$ where
$A=\langle\alpha\rangle$ and $B=\langle\beta\rangle$.
It follows from Eq.~(11) that
$N_A(\lambda_{\alpha,\beta})=N_B(\lambda_{\alpha,\beta})=1$,
and hence a necessary condition is that
$\lambda_{\alpha,\beta}\in K_A\cap K_B$.

\unitlength=1mm
\begin{center}
\begin{picture}(50,85)(0,0)
\put(20,0){\put(-1,-1){$\bullet$ 1}}
\put(20,10){\put(-1,-1){$\bullet$ $I_A\cap I_B=I_{AB}$}}
\put(30,20){\put(-1,-1){$\bullet$ $I_B$}}
\put(10,20){\put(-1,-1){$\bullet$}}
\put(10,20){\put(-6,-1){$I_A$}}
\put(20,30){\put(-1,-1){$\bullet$ $I_A I_B$}}
\put(20,0){\line(0,1){10}}
\put(20,10){\line(1,1){10}}
\put(30,20){\line(-1,1){10}}
\put(20,30){\line(-1,-1){10}}
\put(10,20){\line(1,-1){10}}
\put(20,40){\put(-0.5,-1){$\vdots$}}
\put(20,80){\line(0,-1){10}}
\put(20,70){\line(1,-1){10}}
\put(30,60){\line(-1,-1){10}}
\put(20,50){\line(-1,1){10}}
\put(10,60){\line(1,1){10}}
\put(20,80){\put(-1,-1){$\bullet$ $\Etimes=\Ftimes_{q^{mn}}$}}
\put(20,70){\put(-1,-1){$\bullet$ $K_A K_B=K_{AB}$}}
\put(10,60){\put(-1,-1){$\bullet$}}
\put(10,60){\put(-8,-1){$K_A$}}
\put(30,60){\put(-1,-1){$\bullet$ $K_B$}}
\put(20,50){\put(-1,-1){$\bullet$ $K_A\cap K_B$}}
\end{picture}
\hskip10mm
\begin{picture}(50,85)(0,0)
\put(20,0){\put(-1,-1){$\bullet$ 1}}
\put(20,10){\put(-1,-1){$\bullet$ $q-1$}}
\put(30,20){\put(-1,-1){$\bullet$ $q^n-1$}}
\put(10,20){\put(-1,-1){$\bullet$}}
\put(10,20){\put(-14,-1){$q^m-1$}}
\put(20,30){\put(-1,-1){$\bullet$ $\ell$}}
\put(20,0){\line(0,1){10}}
\put(20,10){\line(1,1){10}}
\put(30,20){\line(-1,1){10}}
\put(20,30){\line(-1,-1){10}}
\put(10,20){\line(1,-1){10}}
\put(20,40){\put(-0.5,-1){$\vdots$}}
\put(20,80){\line(0,-1){10}}
\put(20,70){\line(1,-1){10}}
\put(30,60){\line(-1,-1){10}}
\put(20,50){\line(-1,1){10}}
\put(10,60){\line(1,1){10}}
\put(20,80){\put(-1,-1){$\bullet$ $e=q^{mn}-1$}}
\put(20,70){\put(-1,-1){$\bullet$ $e/(q-1)$}}
\put(10,60){\put(-1,-1){$\bullet$}}
\put(10,60){\put(-20,-1){$e/(q^n-1)$}}
\put(30,60){\put(-1,-1){$\bullet$ $e/(q^m-1)$}}
\put(20,50){\put(-1,-1){$\bullet$ $e/\ell$}}
\end{picture}
\end{center}
\vskip2mm
\begin{center}
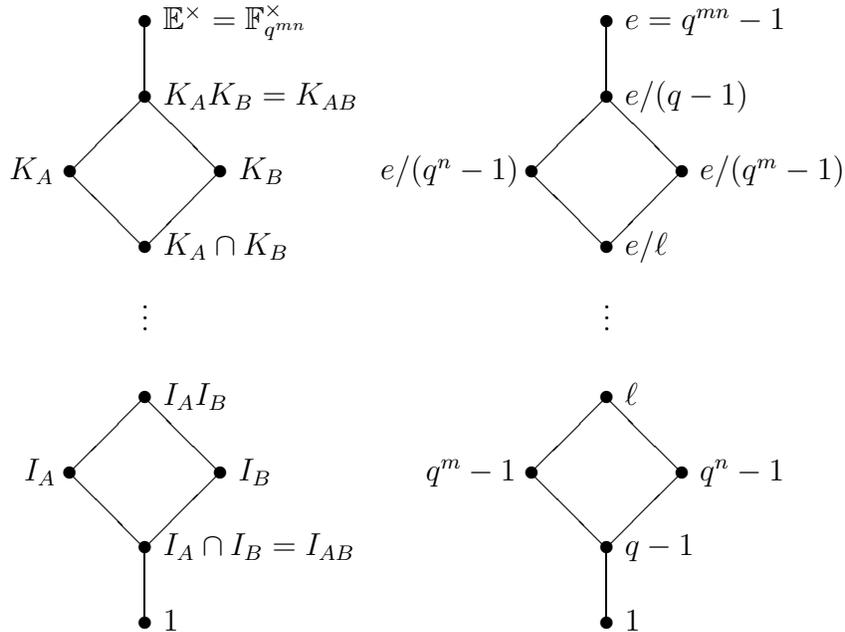
Figure 1. Subgroups of $\Ftimes_{q^{mn}}$ and their orders.
\end{center}
\vskip2mm

Since $AB=BA$ and $A\cap B=1$, $N_{AB}$ equals $N_A\circ
N_B=N_B\circ N_A$, and hence $K_A K_B\subseteq K_{AB}$ and
$I_{AB}\subseteq I_A\cap I_B$. If $\E$ is finite, then these containments
are equalities, and the necessary condition
$\lambda_{\alpha,\beta}\in K_A\cap K_B$ is sufficient to solve
Eq.~(11) for $\nu_\alpha\in K_A$ and $\nu_\beta\in K_B$ (see Section~6).

When $\F=\E^{AB}$ is finite of order~$q$, then $\E=\F_{q^{mn}}$, and
$G=AB$. In Figure~1, $\ell=\gcd(q^m-1,q^n-1)=(q^m-1)(q^n-1)/(q-1)$ and
$I_{AB}=\Ftimes$.

Suppose that $\E$ is a number field and $\rho$ maps into
$\GL_d(\Z_\E)$, where $\Z_\E$ denotes the ring of integers of $\E$.
Then there exist algorithms \cite{C93} for computing the group
$U(\Z_\E)$ of units of $\Z_\E$. Therefore solving
\[
  \beta(\nu_\alpha)^{-1}\nu_\beta^{-1}\nu_\alpha\alpha(\nu_\beta)
  =\lambda_{\alpha,\beta}\qquad(\alpha,\beta\in S)
\]
for $\nu_\alpha\in K_{\langle\alpha\rangle}$,
$\nu_\beta\in K_{\langle\beta\rangle}$
reduces to solving a linear system over $\Z$.

Although evaluating $\Pi_C(X)$ is clearly useful, it is
time-consuming when $|G|$ is large unless an averaging argument is
used. We describe how to use a subgroup chain
$G=G_0>G_1>\cdots>G_{t+1}=1$ to reduce the cost of computing $\Pi_C(X)$ from
$\textup{O}(|G|)$ to $\textup{O}(\sum_{i=1}^t|G_{i-1}: G_i|)$.
If $G=\alpha_1 H\cup\cdots\cup\alpha_r H$ is a
decomposition of $G$ into left cosets of $H$, then
\[
 \Pi_C(X)
 =\sum_{i=1}^r\sum_{\beta\in H}C_{\alpha_i\beta} \alpha_i\beta(X)
 =\sum_{i=1}^r C_{\alpha_i}\alpha_i
    \left(\sum_{\beta\in H}C_{\beta} \beta(X)\right).
\]
Put differently, $\Pi_{C|G}=\sum_{i=1}^r C_{\alpha_i}\alpha_i
\Pi_{C|H}$. If $G$ is solvable, then we may choose $G_i$ so that
$G_i\triangleright G_{i+1}$ and
$G_i=\langle\gamma_i,G_{i+1}\rangle$. In this case, an idea in
\cite[p.\;1705]{GH97} further reduces the complexity of evaluating
$\Pi_C(X)$ to $\textup{O}(\log|G|)$.

\section{Las Vegas algorithms}

A Las Vegas algorithm is one that involves random choices, and when
it terminates it produces an answer that is provably correct.  For
example, a 1-cocycle $C\colon G\to\GLdE$ may be written as
$C_\alpha=A\alpha(A)^{-1}$, $\alpha\in G$, by repeated selecting a
random $X\in\MatdE$ until $A=\Pi_C(X)$ is invertible. If $|\F|=q$ is
finite and a uniform distribution is used for $\MatdE$, then the
probability that $\Pi_C(X)$ is invertible is
\[
 f(d,q)=\frac{|\GLdF|}{|\MatdF|}=\prod_{i=1}^d (1-q^{-i}).
\]
Note that
\[
 \limsup_q f(d,q)=f(d,\infty)=1\quad\text{and}\quad
 \liminf_{d,q} f(d,q)=f(\infty,2).
\]
The following bounds for $f(d,q)$ are useful:
\[
 1-q^{-1}\ge f(d,q)>\prod_{i=1}^\infty (1-q^{-i})
 >1-\sum_{i=1}^\infty q^{-i}=1-(q-1)^{-1}.
\]
If $q=2$, then $f(\infty,2)=0.288788\cdots>2/7$ gives a better lower
bound. Thus one would expect to make on average at most 3.5 choices
for $X$. The probability that the algorithm fails to terminate after
$n$ selections is $(1-f(d,q))^n<\min\{(q-1)^{-n},(5/7)^n\}$.
If $\E$ is infinite, then it follows by localization and a
local-global argument that the
probability that $\Pi_C(X)$ is invertible is 1.

In the light of Theorem~10 we should also consider the probability,
$p_C$, that a random $\lambda\in\Etimes$ has $\Pi_C(\lambda I)$
invertible. If $\E$ is finite, then certain choices for $C$ have
$p_C=1$. Empirical evidence suggests that when $|\E|$ is small
the average value of $p_C$,
averaged over all 1-cocycles $C$, is a number very close to
$f(d,q)$. This is our default expectation.

We describe a Las Vegas algorithm for computing $(q-1)$th roots.
Let $C\colon G\to\textup{GL}_1(\E)$ be a 1-cocycle where
$\E=\F_{q^n}$, $\F=\F_q$ and $G=\langle\alpha\rangle$ where
$\alpha(\lambda)=\lambda^q$. If $C_\alpha=\lambda$, then
$\lambda\alpha(\lambda)\cdots\alpha^{n-1}(\lambda)=1$ and finding
$\mu\in\Etimes$ such that $\lambda=\mu\alpha(\mu)^{-1}$ is equivalent
to finding a $(q-1)$th root, as $\mu^{q-1}=\lambda^{-1}$. Lemma~2(b)
gives a Las
Vegas algorithm for computing $\mu$: choose $\nu\in\E$ randomly
until $\Pi_C(\nu)$ is nonzero. As $\Pi_C$ is a nonzero $\F$-linear
map $\E\to\F$, each $\nu$ has probability $1-q^{-1}$ of success.
Unless $q$ is small, this Las Vegas algorithm is faster than factoring
the polynomial $x^{q-1}-\lambda^{-1}$ over $\E$.

We comment now on Las Vegas algorithms for solving norm equations
in finite fields. Let $\E=\F_{q^n}$, $\F=\F_q$ and let
$\lambda\in\Ftimes$. Denote by $e$, $f$ and $|\lambda|$ the orders of
$\Etimes$, $\Ftimes$ and $\langle\lambda\rangle$ respectively.
One may solve the equation $N(\mu)=\lambda$ by randomly selecting
$\nu\in\Etimes$ and checking whether or not $\mu=\nu^{f/|\lambda|}$
satisfies $N(\mu)=\lambda$. As the norm map
$N\colon\Etimes\to\Ftimes\colon\mu\mapsto\mu^{e/f}$ is surjective,
each selection has probability $|\lambda|^{-1}$ of success.
This algorithm is useful when $|\lambda|$ is small. If $|\lambda|$ is
large, then another Las Vegas algorithm is more desirable. Let
$d=\gcd(|\lambda|,e/f)$. Since $q\equiv 1\mod |\lambda|$, it follows
that $e/f\equiv n\mod |\lambda|$. In most applications, $n$ is small
when $f$ is large, and hence when $|\lambda|$ is large $d$ is usually
much smaller. Denote by $s$ a multiplicative inverse of $e/(fd)$ modulo
$|\lambda|/d$. Randomly select $\nu\in\Etimes$. A root $\mu$ of the
polynomial $x^d-\lambda^s\nu^f$ has probability $d^{-1}$ of satisfying
$N(\mu)=\lambda$.

We prove that the above algorithm is correct, and that either all $d$th
roots $\mu$ of $\lambda^s\nu^f$ satisfy $N(\mu)=\lambda$, or none do.
Let $\Etimes=\langle\zeta\rangle$, and suppose that
$\lambda=\zeta^{ie/f}$ and $\nu=\zeta^j$. As $\zeta^{e/f}$ has order
$f$, $|\lambda|$ equals $f/\gcd(i,f)$. Let $r,s\in\Z$ satisfy
$r|\lambda|+se/f=d$ where $d=\gcd(|\lambda|,e/f)$. If $\mu=\zeta^k$,
then modulo $e$
\begin{align*}
 dk&\equiv ise/f+jf\\
   &\equiv i(d-r|\lambda|)+jf\\
   &\equiv id + (-t+j)f\quad\textup{where }t:=ir/\gcd(i,f)\in\Z.
\end{align*}
There exists an $\ell\in\Z$ such that
\begin{align*}
  k&=i+ (j-t)f/d+\ell e/d\\
 ke/f&= ie/f+(j-t)e/d+\ell (e/d)(e/f)\\
     &\equiv ie/f+(j-t)e/d\mod e.
\end{align*}
So
$N(\mu)=N(\zeta^k)=\zeta^{ke/f}=\zeta^{ie/f}(\zeta^{e/d})^{j-t}
=\lambda\omega^{j-t}$
where $\omega=\zeta^{e/d}$ has order $d$. In summary, $N(\mu)=\lambda$
if and only if $j\equiv t\mod d$. Thus the probability of success is
$d^{-1}$. As the value of $N(\mu)$ is independent of $\ell$, either
each of the $d$ roots $\mu$ satisfy $N(\mu)=\lambda$, or none do.

In the case when $\E$ is finite and $|G|=|\E:\F|$ is not a prime
power, then a divide-and-conquer strategy may be used for solving norm
equations. Suppose that $|\E:\F|=mn$ where $\gcd(m,n)=1$ and
$\F=\F_q$. Let $G=AB$ where $A=\langle\alpha\rangle$ satisfies
$\alpha(\lambda)=\lambda^{q^n}$, and $B=\langle\beta\rangle$
satisfies $\beta(\lambda)=\lambda^{q^m}$. Then $|A|=m$ and $|B|=n$.
If the presentation
$G=\langle\alpha\beta\mid(\alpha\beta)^{mn}=1\rangle$ is used,
then one need only solve one norm equation:
$\mu_{\alpha\beta}^{(q^{mn}-1)/(q-1)}=\lambda_{\alpha\beta}$ where
$\lambda_{\alpha\beta}\in\Ftimes$ is given. If the
presentation
$G=\langle\alpha,\beta\mid\alpha^m=\beta^n=[\beta,\alpha]=1\rangle$ is
used, then one must solve three equations:
$\mu_{\alpha}^{(q^{mn}-1)/(q^n-1)}=\lambda_\alpha$,
$\mu_{\beta}^{(q^{mn}-1)/(q^m-1)}=\lambda_\beta$ and
$\nu_\alpha^{1-q^m}\nu_\beta^{q^n-1}=\lambda_{\alpha,\beta}$ where
$\lambda_\alpha\in\E^A$, $\lambda_\beta\in\E^B$ and
$\lambda_{\alpha,\beta}\in K_A\cap K_B$. The two norm equations could
be solved using the above Las Vegas algorithm. This has the advantage
that $\gcd(|\lambda_\alpha|,m)$ and $\gcd(|\lambda_\beta|,n)$ are
likely smaller than $\gcd(|\lambda_{\alpha\beta}|,mn)$. There exist
$r,s\in\Z$ such that
\[
 r(q^m-1)+s(q^n-1)=q-1.
\]
Since $\lambda_{\alpha,\beta}\in K_A\cap K_B\subseteq K_{AB}$, our Las
Vegas algorithm for computing $(q-1)$th roots may be used to solve the
equations $\nu_\alpha^{q-1}=\lambda_{\alpha,\beta}^{-r}$ and
$\nu_\beta^{q-1}=\lambda_{\alpha,\beta}^s$. Then
\[
 \nu_\alpha^{1-q^m}\nu_\beta^{q^n-1}=
 \lambda_{\alpha,\beta}^{r(q^m-1)/(q-1)+s(q^n-1)/(q-1)}
 =\lambda_{\alpha,\beta}.
\]

\section{Remarks and examples}

The assumption that $\rho$ is absolutely irreducible was not used in
Sections 1--4, however, it is very useful for practical algorithms for
writing $\rho$ over $\F$. If $\rho$ is reducible, then one may need to
solve linear systems to find $D$ satisfying Eq.~(9), and the solution spaces
may be more than one-dimensional. Finding $C$ from
$D$ is likely to be problematic. If $\rho$ is irreducible but
not absolutely irreducible, then the \Meataxe\ \cite{HR94,P98} may be
used to find $D$. In this case, however, the arithmetic
needed to solve for $\mu$ (and hence find $C$)
takes place in the division algebra $\End(\rho)$ of matrices commuting
with $\rho(\A)$. See \cite{G03} for a description of some of the
relevant noncommutative theory. We shall assume henceforth that
$\A=\FH$ is a group algebra.

The connection between $\EH$-modules and $\FH$-modules is clarified by
considering normal bases. The following simple observation is not made
explicitly in texts covering modular representation theory such as
\cite{HB82}. Let $(\alpha(\lambda))_{\alpha\in G}$ be a normal basis
for $\E$ over
$\F$. Let $V=\E^{d\times 1}$ and $U=\F^{d\times 1}$. Then $V$ viewed
as an $\FH$-module is a direct sum of $|G|$ Galois conjugate
$\FH$-submodules: $V=\dotplus_{\alpha\in G} \alpha(\lambda) U$.
Note that $A^{-1}\rho(h)A\in\GLdF$ for $h\in H$ and so
\[
 \alpha(\lambda) UA^{-1}\rho(h)A=\alpha(\lambda) U=\alpha(\lambda U).
\]
Thus the $\alpha(\lambda) U$ are $A^{-1}\rho A$ invariant, and Galois
conjugate.

In the examples below $\E=\F(\zeta_n)$ is a subfield of the complex
numbers, and $\zeta_n=e^{2\pi i/n}$. An automorphism
$\alpha$ of $\E$ is determined by a number $k$ satisfying
$\alpha(\zeta_n)=\zeta_n^k$ and $\gcd(k,n)=1$. As usual, $\Q$
denotes the rational field.

\subsection*{Example 1} Let $H$ be the dicyclic group of order $8n$
\[
 H=\langle a,b\mid a^2=b^{2n}, b^{4n}=1,a^{-1}ba=b^{-1}\rangle.
\]
Let $\E=\Q(\zeta)$ where $\zeta=\zeta_{4n}$. Define
$\alpha\in\Aut(\E)$ by $\alpha(\zeta)=\zeta^{-1}$. Then
$\alpha$ has order~2, and
$\F=\E^{\langle\alpha\rangle}=\Q(\zeta+\zeta^{-1})$. Define
$\rho\colon H\to\GL_2(\E)$ by
\[
 \rho(a)=\begin{pmatrix}0&1\\-1&0\end{pmatrix}\quad\text{and}\quad
 \rho(b)=\begin{pmatrix}\zeta&0\\0&\zeta^{-1}\end{pmatrix}.
\]
Then $D_\alpha=\rho(a)$ and $D_1=\rho(1)$ satisfies Eq.~(9).
Since $N_{\langle\alpha\rangle}(D_\alpha)$ equals
$D_\alpha\alpha(D_\alpha)=D_\alpha^2=-I$, it follows that
$\lambda_\alpha=-1$. Since $\alpha$ is complex conjugation,
$N_{\langle\alpha\rangle}(\mu_\alpha)=\mu_\alpha\overline{\mu_\alpha}
=||\mu_\alpha||^2\ge0$, so $N_{\langle\alpha\rangle}(\mu_\alpha)=-1$ has no
solution. Consequently, $\rho$ can not be written over $\F$.

\subsection*{Example 2} Let $H=\langle a,b\mid a^2=b^{4n}, b^{8n}=1,
a^{-1}ba=b^{1+4n}\rangle$ and let $\E=\Q(\zeta)$ where $\zeta=\zeta_{8n}$.
Define $\alpha\in\Aut(\E)$ by $\alpha(\zeta)=\zeta^{1+4n}=-\zeta$.
Then $\alpha$ has order~2, and
$\F=\E^{\langle\alpha\rangle}=\Q(\zeta^2)$. Define
$\rho\colon H\to\GL_2(\E)$ by
\[
 \rho(a)=\begin{pmatrix}0&1\\-1&0\end{pmatrix}\quad\text{and}\quad
 \rho(b)=\begin{pmatrix}\zeta&0\\0&\zeta^{1+4n}\end{pmatrix}.
\]
Set $D_1=\rho(1)$ and $D_\alpha=\rho(a)$. Then
$N_{\langle\alpha\rangle}(D_\alpha)=-I$, so $\lambda_\alpha=-1$.
Now $\mu_\alpha=\zeta^{2n}$ satisfies
$N_{\langle\alpha\rangle}(\mu_\alpha)=\mu_\alpha^2=-1=\lambda_\alpha^{-1}$.
Thus $C_1=\rho(1)$ and $C_\alpha=\zeta^{2n}\rho(a)$. The matrix
\[
   A:=\Pi_C\left(\frac{1+\zeta}2\,I\right)
    =\frac{1}2\begin{pmatrix}1+\zeta&\zeta^{2n}(1-\zeta)\\
         -\zeta^{2n}(1-\zeta)&1+\zeta\end{pmatrix}
\]
has $\det(A)=\zeta\ne0$, and hence writes $\rho$ over $\F$. If
$\rho'=A^{-1}\rho A$, then  
\[
 \rho'(a)=\begin{pmatrix}0&1\\-1&0\end{pmatrix}\quad\text{and}\quad
 \rho'(b)=\frac{1}2\begin{pmatrix}1+\zeta^2&\zeta^{2n}(1-\zeta^2)\\
         \zeta^{2n}(1-\zeta^2)&-1-\zeta^2\end{pmatrix}.
\]
The similarity between $A$ and $\rho'(b)$ is interesting. For each $n$
there are many choices for $\mu_\alpha$, and then many choices for
$\nu$ such that $\Pi_C(\nu I)$ is invertible. Our choices
$\mu_\alpha=\zeta^{2n}$, $\nu=(1+\zeta)/2$ give a simple expression
for $\rho'(b)$. Another choice when $n$ is odd is 
$\mu_\alpha=1+\zeta^n-\zeta^{3n}$ and $\nu=1$.

\subsection*{Example 3} Let $H=\langle a,b\mid a^m=b^n=1, a^{-1}ba=b^r\rangle$
where $r$ is the order of $m$ modulo $n$. Let $\zeta=\zeta_n$,
$\E=\Q(\zeta)$, and let $\F=\E^{\langle\alpha\rangle}$ where
$\alpha\in\Aut(\E)$ is defined by $\alpha(\zeta)=\zeta^r$. Define
$\rho\colon H\to\GL_m(\E)$ by
\[
 \rho(a)=
 \begin{pmatrix}0&1 &0\\ & &\ddots&\\0&0& &1\\1&0& &0\end{pmatrix}
 \quad\text{and}\quad
 \rho(b)=\begin{pmatrix}\zeta& & & \\ &\zeta^r& & \\ & &\ddots
  &\\&&&\zeta^{r^{m-1}}\end{pmatrix}.
\]
Then $C_\alpha=\rho(a)$ and
$C_{\alpha^i}=C_\alpha\alpha(C_\alpha)\cdots\alpha^{i-1}(C_\alpha)=\rho(a)^i$
and
\[
 A=\Pi_C(\lambda I)=\sum_{i=0}^{m-1} C_\alpha^i\alpha^i(\lambda)
 =(\alpha^{i-j}(\lambda))
\]
is invertible if and only if $\lambda$ defines a
normal basis for $\E$ over $\F$.
If $\rho'=A^{-1}\rho A$, then $\rho'(a)=\rho(a)$ and the expression
for $\rho'(b)$ is rather complicated, and depends on $r$.

\subsection*{Example 4} Let $\E:\F$ be a finite Galois extension with
group $G$. Let $\sigma$ be the left regular representation
$G\to\textup{Sym}(G)$ satisfying
$\sigma_\alpha(\gamma)=\alpha\gamma$ and
$\sigma_{\alpha\beta}=\sigma_\alpha\circ\sigma_\beta$. Let $H$ be the split
extension of $\Etimes$ by $G$. Specifically, let $H=G\times\Etimes$
where
\[
 (\alpha,\lambda)(\beta,\mu)=(\alpha\beta,\beta(\lambda)\mu)\qquad
 (\alpha,\beta\in G,\lambda,\mu\in\Etimes).
\]
Define $\rho\colon H\to\GL_{|G|}(\E)$ by
$\rho(\alpha,\lambda)=(\eta(\lambda)\delta_{\sigma_\alpha(\eta),\eta})$
where $(\delta_{\xi,\eta})$ is the identity matrix. The $(\xi,\eta)$
entry of $\rho(\alpha,\lambda)$ is zero unless
$\xi=\sigma_\alpha(\eta)$ in which case it equals $\eta(\lambda)$.
The $(\xi,\eta)$ entry of $\rho(\alpha,\lambda)\rho(\beta,\mu)$ is zero unless
$\xi=\sigma_{\alpha\beta}(\eta)$ in which case it equals
$\sigma_\beta(\eta)(\lambda)\eta(\mu)=\eta(\beta(\lambda)\mu)$.
This proves that $\rho$ is a homomorphism. Since $\rho$ is induced
from a 1-dimensional representation $\Etimes\to\GL_1(\E)$ which is
fixed only by the identity automorphism, it follows from Clifford's
theorem that $\rho$ is absolutely irreducible. We may take
$C_\alpha$ to be the permutation matrix $\rho(\alpha,1)$ corresponding
to $\sigma_\alpha$. Then $A=\Pi_C(\lambda I)$ is invertible if and only if
$\lambda$ defines a normal basis for $\E$ over $\F$. If $|\F|=q$ and
$|\E|=q^n$, then the probability that $\Pi_C(\lambda I)$ is invertible
is $q^{-n}\sum_{d|n}\mu(n/d)q^d$ where $\mu$ denotes the
M\"obius function. It follows by considering base-$q$ expansions that
$\sum_{d|n}\mu(n/d)q^d\ge q^n-q^{n/p}\ge q^n-q^{n/2}$
where $p$ is the smallest prime divisor of $n$. Hence
$q^{-n}\sum_{d|n}\mu(n/d)q^d\ge 1-q^{-n/2}$.

\frenchspacing
{\renewcommand{\baselinestretch}{1}

} 

\vskip3mm
\goodbreak
\def\efont{\scriptsize\upshape\ttfamily}
{\tiny\scshape
\begin{tabbing}
\=\hspace{70mm}\=\kill\\
\>Department of Mathematics    \\
\>Central Washington University\\
\>WA 98926-7424, USA           \\
\>\efont GlasbyS@cwu.edu       \\
\end{tabbing}
}

\end{document}